\documentclass{article}
\usepackage{amsmath,amssymb}

\newcommand{\F}{\mathcal{F}}
\newcommand{\Fi}{\overline{F}_\infty}

\newcommand{\Fen}{\overline{F}_{\eps_n}}
\newcommand{\Fn}{\overline{F}_{n}}
\newcommand{\tF}{\tilde{F}}

\newcommand{\dis}{\displaystyle}

\newcommand{\R}{\mathbb{R}}

\newcommand{\C}{\mathcal{C}}

\newcommand{\M}{{\mathcal{M}}}
\newcommand{\Ou}{{\mathcal{O}}}
\newcommand{\D}{{\mathcal{D}}}
\newcommand{\sign}{{\mathrm{Sign}}}
\newcommand{\eps}{\varepsilon}
\newcommand{\Fe}{{\overline{F}_\varepsilon}}
\newcommand{\fe}{{f_\varepsilon}}
\newcommand{\dt}{\partial_t}
\newcommand{\ds}{\partial_s}
\newcommand{\dy}{\partial_y}
\newcommand{\dx}{\partial_x}

\newcommand{\Lloc}{L_{\mathrm{loc}}}
\def\qed{\hbox{${\vcenter{\vbox{
  \hrule height 0.4pt\hbox{\vrule width 0.4pt height 6pt
  \kern5pt\vrule width 0.4pt}\hrule height 0.4pt}}}$}}

\newtheorem{theo}{Theorem}
\newtheorem{prop}[theo]{Proposition}
\newtheorem{lemm}[theo]{Lemma}

\author{Alexis Vasseur \thanks{Department of Mathematics, University of Texas}}
\title{Rigorous derivation of the Kinetic/Fluid coupling involving a Kinetic layer on
          a toy problem}
\date{}
\begin{document}
\maketitle
\bibliographystyle{plain}
\noindent{\bf Abstract:} In this article, we investigate the
kinetic/fluid coupling on a toy model. We obtain it rigorously from
a hydrodynamical limit. The idea is that at the level of the full
kinetic model, the coupling is obvious. We then investigate the
coupling obtained when passing into the limit. We show that,
especially  in presence of a shock stuck  on the interface, the
coupling involves a kinetic layer known as the Milne problem. Due to
this layer, the limit process is quite delicate and some blow-up
techniques are needed to ensure its strong convergence. \vskip0.3cm
\noindent {\bf MSC-class: 82C40, 35L65, 82D05} 

\section{Introduction}

Fluid models of mechanics, like the compressible Euler system, are
known to be not valid for gas far from thermodynamical equilibrium.
This is especially the case for rarefied gas, or gas in the context
of Physics of high energy. In this case, it is necessary to consider
the kinetic theory which involves a new 3D variable called variable
of velocity, e.g.,  the Boltzmann equation. In a
numerical point of view, this makes the system far more costly to
compute. Different strategies have been developed  to compute with
kinetic theory only in a small region where the gas is known to be
very far from equilibrium. Standard fluid schemes are then used in
the rest of the space. We then have the question of the coupling
between those two regions of the two kind of models.

These kind of coupling conditions have been widely studied,
especially for aerospace engineering. This is the case, in
particular, for the coupling of the Boltzmann equation and
Navier-Stokes system (see, e.g.,  J.-F.Bourgat, P.Le Tallec,
B.Perthame and Y.Qiu \cite{N1},  P.Le Tallec et F.Mallinger
\cite{N2}, or S.Dellacherie \cite{SD}). Most of the methods use the
thermo-equilibrium kinetic function (the Maxellian function) to
compute the in-going fluxes at the boundary.

On the other hand, in the case of the evaporation/condensation  problem
 it is well known that a layer problem is involved. In
the case of Boltzmann equation, this layer problem is called Milne
problem. This problem has been extensively studied by the team of
Y.Sone and K.Aoki (see \cite{SA}). A.V. Bobylev, R.
Grzhibovskis, and A. Heintz proposed a good approximation of it in
\cite{Boby}. Solutions of this kind of equation can be classified in
different categories depending on the boundary conditions. The
complete classification was obtained by  F.Golse, F.Coron and
C.Sulem in \cite{GS}. Finally, study of such kinetic layer can be
found in asymptotic problem linked to Semi-conductors (see
 P.Degond and C.Schmeiser \cite{SC1}, or  N.Ben Abdallah,
P.Degond and I.M.Gamba \cite{SC2}).

In this paper we consider a simplified model for which we will be
able to obtain rigorously  the ``good" coupling by performing an
hydrodynamical limit in one region. The model was introduced by
B.Perthame and E.Tadmor \cite{PT}. It is a caricature of the
Boltzmann equation and is given by:
\begin{equation}\label{BGK}
\dt f+\xi\dx f=\M f-f,
\end{equation}
where the ``equilibrium function" $\M f$ is defined by:
\begin{equation}\label{maxwellienne}
\begin{array}{l}
\M f(t,x,\xi)=\dis{{\bf 1}_{\{0\leq \xi\leq u(t,x)\}}{\bf
-1}_{\{-u(t,x)\leq \xi\leq0 \}}} \\[3mm]
\mathrm{with:} \ \
\dis{u(t,x)=\int_{-\infty}^{+\infty}f(t,x,\xi)\,d\xi}.
\end{array}
\end{equation}
The Boltzmann collision operator is replaced by the relaxation term
$\M f-f$. This equation is a ``BGK" version of the ``transport-collapse" method introduced independently by Y.Brenier
\cite{B1,B} and Y.Giga and Y.Miyakawa \cite{GM} as a numerical
scheme for scalar conservation laws.

Consider a function $u^0\in L^\infty(\R)$. We can define the
associated equilibrium function $f^0=M(u^0(x);\xi)\in
L^\infty(\R\times\R)$ such that:
\begin{eqnarray*}
&&\M f^0=f^0\\
&&\int_\R f^0(x,\xi)\,d\xi=u^0(x) \ \ \mathrm{for  \ every \ } x\in
\R.
\end{eqnarray*}
In  \cite{LPT}, P.L. Lions, B.Perthame and E.Tadmor showed that the
solutions  $f_\eps\in L^\infty(\R^+\times\R\times\R)$ of the
rescaled  equation (\ref{BGK})
\begin{equation}\label{BGKrescale}
\dt \fe+\xi\dx \fe=\frac{\M \fe-\fe}{\eps},
\end{equation}
with initial value $M(u^0;\cdot)$ converge, when $\eps$ goes to 0,
to $M(u(t,x);\xi)$, where $u\in L^\infty(\R^+\times\R)$ is solution
to the Burgers equation
\begin{equation}\label{burgers}
\begin{array}{l}
\dis{\dt u+\dx\frac{u^2}{2}=0,} \\[3mm]
u(t=0)=u^0,
\end{array}
\end{equation}
and verifies the entropy conditions
\begin{equation}\label{entropy}
\dt\phi(u)+\dx H(u)\leq0,
\end{equation}
for every convex function $\phi$ with associated entropy flux $H$
verifying $H'(y)=y\phi'(y)$. This equation mimics the Euler equation
of the compressible gas. In this article, our aim is to derive a
natural coupling condition between equations (\ref{BGK}) and
(\ref{burgers}). We choose to obtain it from the asymptotic limit
when $\eps$ converges to 0 of the following problem
\begin{equation}\label{couplagecinetique}
\begin{array}{l}
\dis{\dt \fe+\xi\dx \fe=\alpha_\eps(x)(\M \fe-\fe), \ \ \mathrm{for} \ t\in\R^+, x\in\R,
\xi\in\R,}\\[0.3cm]
\dis{\fe(t=0)=f^0, \ \ \mathrm{for} \ x\in\R, \xi\in\R,}
\end{array}
\end{equation}
where the function $\alpha_{\eps}(x)$ is $1$ for $x\leq0$ and
$1/\eps$ for $x>0$. This corresponds to solving the problem
(\ref{BGK}) in the domain $x<0$ and the problem (\ref{BGKrescale})
in the domain $x>0$ with a boundary coupling of the type complete
transmission at $x=0$.

In most of the numerical coupling, the in-flux entering in the
kinetic part is defined through the equilibrium function. The
coupling, in this simplified model would be
\begin{equation}\label{couplagenumerique}
\begin{array}{l}
\dis{\F^+(u(t,0+))=\int_0^\infty \xi f(t,0-,\xi)\,d\xi}\\[0.3cm]
\dis{f(t,0-,\xi)=M(u(t,0+);\xi) \ \ \mathrm{for} \ \xi\leq0,}
\end{array}
\end{equation}
where  $\F^+$ is the positive semi-flux of the Engquist-Osher scheme
\cite{EO} (see also Degond and Jin \cite{DegondShi} for another
strategy of coupling involving a smooth transition). This provides a
plausible coupling for the the system (\ref{BGK}) (\ref{burgers}).
However, this is not the one obtained in the hydrodynamical limit.
This is due to the possible production of a kinetic layer of size
$\eps$ at the interface $x=0$ (the equivalent of the one described
in the context of Boltzmann equation by the Milne problem)  (see
A.Klar \cite{N3}). We show, in this article, that the correct
coupling obtained at the limit of (\ref{couplagecinetique}) is the
following one
\begin{equation}\label{limitefluide}
\left\{
\begin{array}{l}
\dis{\dt u+\dx\frac{u^2}{2}=0,} \ \ \mathrm{for} \
t\in\R^+,x\in\R^+, \\[3mm]
u(t=0,x)=u^0(x)=\int_{-\infty}^{+\infty}f^0(x,\xi)\,d\xi, \ \
\mathrm{for} \ x\in\R^+,
\\[3mm]
\dis{u(t,x=0+)\stackrel{BLN}{:=}\sqrt{2\int_0^\infty\xi
f(t,0-,\xi)\,d\xi}, \ \ \mathrm{for} \ \ t\in\R^+,}
\end{array}
\right.
\end{equation}
\vskip2mm
\begin{equation}\label{limitecinetique}
\left\{
\begin{array}{l}
\dt f+\xi\dx f=\M f-f, \ \ \mathrm{for} \
t\in\R^+,x\in\R^-,\xi\in\R,\\[3mm]
f(t=0,x,\xi)=f^0(x,\xi)  \ \ \mathrm{for} \ x\in\R^-,\xi\in\R,\\[3mm]
f(t,x=0-,\xi\leq0)=F(t,y=0+,\xi\leq0) \ \ \mathrm{for} \
t\in\R^+,\xi\in\R^-,
\end{array}
\right.
\end{equation}
\vskip2mm
\begin{equation}\label{limitecouche}
\left\{
\begin{array}{l}
\xi\dy F=\M F-F, \ \ \mathrm{for} \
t\in\R^+,y\in\R^+,\xi\in\R,\\[3mm]
F(t,y=0+,\xi\geq0)=f(t,x=0-,\xi\geq0), \ \ \mathrm{for} \
t\in\R^+,\xi\in\R^+,\\[3mm]
\dis{\int_{-\infty}^{+\infty}\xi
F(t,y,\xi)\,d\xi=\frac{u(t,0+)^2}{2}} \ \ \mathrm{for} \
t\in\R^+,y\in\R^+.
\end{array}
\right.
\end{equation}
Note that the limit system should verify the equality of the flux
at the interface:
\[
\int_{-\infty}^{+\infty}\xi f(t,x=0-,\xi)\,d\xi=\frac{u(t,0+)^2}{2}.
\]
The coupling condition in  (\ref{limitefluide}) shows that it is
possible, (as suggested in  (\ref{couplagenumerique})), to define
the boundary condition for the fluid equation only from the
out-going flux of the kinetic part. It is worth noticing that this
condition corresponds exactly to the condition of
Bardos-Leroux-N\'ed\'elec \cite{BLN}, which is the natural boundary
condition of Dirichlet type for the Burgers equation.
 It is defined in the following
way
\begin{eqnarray*}
&&\qquad\qquad u\stackrel{BLN}{:=}v\\
&& \qquad\qquad\mathrm {if \ \ and \ \ only \ \ if}\\
&& [u^2/2-k^2/2]\sign(u-v)\leq0 \qquad \mathrm{ for \ \ every \ \ }
k \mathrm{\ between \ } u \mathrm{\ and \ } v.
\end{eqnarray*}

This result, in some sense, justifies the relevance of the Milne
problem in the case of evaporation/condensation problems. At the
numerical point of view, the result could seem a little bit odd. The
computation of the Milne problem is very costly, and is certainly
not needed anyway. Indeed, we   show that there exists two types of
layers depending on the boundary values on each side of the
interface. In the case of the rarefaction layer, we show that the
boundary condition coincides with the numerical condition
(\ref{couplagenumerique}). The layer problem provides a different
boundary condition only in the other case of shock layer.  This case
corresponds physically to a shock coming from the fluid domain that
artificially sticks to the interface. In some way, this means that
the position of the interface has not been wisely chosen. To get an
accurate description of the physics, the full shock should be
modeled by the kinetic equation. More sophisticated numerical schemes
use a moving interface chosen on the fly, tracking exactly the
motion of the shocks near this interface (see  Goldstein
and al \cite{Goldstein}).

For this simplified model, first steps were done  by
A.Nouri, A.Omrane and  J.P.Villa \cite{NOV}, C.Bourdarias, M.Gisclon
and A.Omrane \cite{BG} and M.Tidriri \cite{T1,T2}. However, none of
those works consider the existence of the kinetic layer. Let us also
mention that
 F.Golse gave a detailed study of the shock profiles of (\ref{BGK})
 in the whole line (see \cite{G}).

We give in the next section the precise statements of the results.
%Note that we do not have a  analysis in full generality of the
%asymptotic limit of (\ref{couplagecinetique}).
The characterization of the coupling condition in
(\ref{limitefluide}) relies on the result of strong traces for
scalar conservation laws \cite{V3} (see Panov \cite{Panov} and
\cite{KV} for new versions).

The convergence in the kinetic layer is obtained thanks to a ``blow-up" method first applied in the hyperbolic context in
\cite{V1,V2} (see also \cite{Ottonew} for more elaborate results using this technique). We present this technique in  section
\ref{sectionlocal}. It makes use of a Liouville's type lemma presented
in  section \ref{sectionliouville}.

Let us fix $L>0$. From now on, we will consider initial values
supported in $\xi$ in $[-L,L]$. A maximum principle implies that for
every  $\eps>0$ the solutions to  (\ref{limitecouche}), (\ref{couplagecinetique}), and
(\ref{limitecinetique})    preserve this
support property (see Perthame \cite{Perthame}). This condition coincides with  solutions $u$ of
(\ref{limitefluide}) that are bounded by $L$.

The results of this paper were announced some time ago in
\cite{XEDP}.

\section{Statement of the results}\label{sectionpres}

We first show  the well-posedness and stability of the limit
problem (\ref{limitefluide}) (\ref{limitecinetique}) (\ref{limitecouche}).
\begin{theo}\label{theoprolimite}
Let $(f^0,u^0)\in L^\infty(\R^-\times[-L,L])\times L^\infty(\R^+)$
be such that  $0\leq \sign(\xi)f^0\leq1$ and $|u^0|\leq L$. Then
there exists a unique  solution $(f,u,F)\in
L^\infty(\R^+\times\R^-\times[-L,L])\times
L^\infty(\R^+\times\R^+)\times L^\infty(\R^+\times\R^+\times[-L,L])$
to (\ref{limitefluide}) (\ref{limitecinetique})(\ref{limitecouche})
verifying $0\leq\sign(\xi)f\leq1$ and $0\leq\sign(\xi)F\leq1$.
Moreover, two such solutions $(f_1,u_1,F_1)$, $(f_2,u_2,F_2)$ verify
for any  $t\in\R+$:
\begin{eqnarray*}
&&\| f_1(t)-f_2(t)\|_{L^1(\R^-\times[-L,L])}+\|u_1(t)-u_2(t)
\|_{L^1(\R^+)}\\
&&\qquad \leq \| f_1^0-f_2^0\|_{L^1(\R^-\times[-L,L])}+\|u_1^0-u_2^0
\|_{L^1(\R^+)}.
\end{eqnarray*}
\end{theo}
Note that we do not claim the $L^1$ stability of the layer $F$ (Indeed, as we will see later, it is not stable).

Then we consider the asymptotic limit. We prove the asymptotic limit
in two different situations. The first set of initial values
considered is the following. Let $0<u^+<L$ and  consider the functions
$\fe^0\in L^\infty(\R\times[-L,L])$ with $0\leq
\sign(\xi)\fe^0\leq1$ verifying
\begin{equation}\label{initialrelax}
\begin{array}{l}
\fe^0(x,\xi)\geq M(u^+;\xi) \ \ \mathrm{for} \ \
x\in\R^-,\xi\in[-L,L],\\
\fe^0(x,\xi)\geq M(-u^+;\xi) \ \ \mathrm{for} \ \
x\in\R^+,\xi\in[-L,L].
\end{array}
\end{equation}
The second set of initial values is the following.
\begin{equation}\label{initialchoc}
\begin{array}{l}
\mathrm{Supp}f^0(x,\cdot)\subset ]-L,u^+-\eta] \ \ \mathrm{for} \ \
x\in\R^-,\\
f^0(x,\xi)=M(-u^+;\xi) \ \ \mathrm{for} \ \ x\in\R^+,\xi\in[-L,L].
\end{array}
\end{equation}
Note that this case is more restrictive, since it is constant in the
fluid region.

We then prove the following result.
\begin{theo}\label{theolimitasymp}
Consider initial values $\fe^0$ verifying
(\ref{initialrelax}) and converging weakly   to a function $f^0\in
L^\infty(\R\times[-L,L])$ in $L^\infty*$. Let  $\fe$ be the solution
to (\ref{couplagecinetique}) with initial value  $\fe^0$. Then the
family $(\fe|_{\{x\leq0\}},u_\eps=\int \fe\,d\xi|_{\{x\geq0\}})$
converges weakly to the unique solution $(f,u)$ to
(\ref{limitefluide}) (\ref{limitecinetique}) (\ref{limitecouche}),
with initial value  $(f^0|_{\{x\leq0\}},u^0=\int
f^0\,d\xi|_{\{x\geq0\}})$.

Consider now an initial value $f^0$ verifying (\ref{initialchoc}).
Let $\fe$ be the solution to (\ref{couplagecinetique}) associated to
this initial value. Let $F_\eps$ be defined by
\begin{equation}\label{Feps}
F_\eps(t,y,\xi)=\fe(t,\eps y,\xi), \ \ \mathrm{for} \
t\in\R^+,y\in\R^+,\xi\in [-L,L].
\end{equation}
Then $(\fe|_{\{x\leq0\}},u_\eps=\int
\fe\,d\xi|_{\{x\geq0\}},F_\eps)$ converges  to $(f,u,F)$ solution to
(\ref{limitefluide}) (\ref{limitecinetique}) (\ref{limitecouche})
with initial value
\begin{eqnarray*}
&&f^0(x,\xi) \ \ \mathrm{for} \ \ x\in\R^-,\xi\in[-L,L],\\
&&u^0(x)=\int_{-L}^{L}f^0(x,\xi)\,d\xi, \ \ \mathrm{for} \ \
x\in\R^+.
\end{eqnarray*}
\end{theo}
Note that most of the difficulties of this limit lie in  obtaining the
boundary condition on the in-flow fluxes in the kinetic domain.
Indeed, the boundary conditions on the fluid part can be obtained, in
a full generality of initial conditions (see Proposition
\ref{theofluide}).

Note that the initial values of type (\ref{initialrelax})  can
converge only weakly. This means that the oscillations on the
initial values cannot propagate in the nonlinear terms of the kinetic and fluid domains.
However, as we will see later, the kinetic layer can be destroyed in
the limit. Indeed, in this case, we do not claim convergence of the
kinetic layer. But the limit functions behave as if the kinetic
layer was present. For the second set of initial values,  things are
more tedious. In this case, the limit solution depends on the
precise structure of the kinetic layer. It is then crucial to show
its convergence. Let us first give the idea why such a layer should
appear. Consider the function $F_\eps$ defined on
$\R^+\times\R^+\times[-L,L]$ by (\ref{Feps}). This function $F_\eps$
verifies
\begin{eqnarray*}
&&\eps\dt F_\eps+\xi\dy F_\eps=\M F_\eps-F_\eps,  \ \
\mathrm{for} \ t\in\R^+,y\in\R^+,\xi\in[-L,L],\\
&&F_\eps(t,y=0,\xi)=\fe(t,x=0,\xi),   \ \ \mathrm{for} \
t\in\R^+,\xi\in[-L,L].
\end{eqnarray*}
Passing to the limit formally, we get
\begin{equation}\label{formel}
\left\{
\begin{array}{l}
\xi\dy F=\M F-F,  \ \
\mathrm{for} \ t\in\R^+,y\in\R^+,\xi\in[-L,L],\\
F(t,y=0,\xi)=f(t,x=0-,\xi),   \ \ \mathrm{for} \
t\in\R^+,\xi\in[-L,L],\\
 F(t,y=+\infty,\xi)=M(u(t,0+);\xi), \ \
\mathrm{for} \ t\in\R^+,\xi\in[-L,L].
\end{array}
\right.
\end{equation}
The last equation is obtained formally by continuity with the fluid
domain. We will show later that, even if the limit problem can
always be formulated in terms of the kinetic layer,  in some cases this asymptotic
limit on the layer can fail.

In order to grasp this fact, we need  to get a refined study of
the kinetic layer (or Milne problem). Consider the set
\[
\C=\left\{ (V,g)\in[0,L^2/2]\times L^1(\R^+); 0\leq g(\xi)\leq 1;
V\geq\int_0^L\xi g(\xi)\,d\xi\right\}.
\]
For every boundary data $(V,g)\in\C$, we say that   $F\in
L^\infty(\R^+\times\R)$ with $0\leq \sign(\xi)F(y,\xi)\leq1$ is
a solution to the kinetic layer problem with given data $(V,g)$ if $F$
verifies
\begin{equation}\label{Milne}
\left\{
\begin{array}{l}
\xi\dy F=\M F-F \ \ \mathrm{for} \  \ y\in\R+,\xi\in(-L,L),\\[3mm]
F(y=0,\xi\geq0)=g \ \ \mathrm{for} \  \ \xi\in(0,L),\\[3mm]
\dis{\int_{-L}^{L} \xi F(y,\xi)\,d\xi=V \ \ \mathrm{for} \ \
y\in\R+.}
\end{array}
\right.
\end{equation}
We show the following proposition which classifies the solution to
(\ref{Milne}).
\begin{prop}\label{theocouche}
For any   $(V,g)\in\C$, there exists a unique solution to
(\ref{Milne}). Let  $u_\infty\in\R^+$ be such that
\[
V=\frac{u_\infty^2}{2}.
\]
Then there are two cases.

\noindent (a) If $V=\int_0^\infty \xi g(\xi)\,d\xi$, we have
$F(y,\xi)=0$ for every  $\xi\leq0$ and
\[
\lim_{y\to\infty}F(y,\cdot)=M(u_\infty,\cdot) \ \ \mathrm{in} \
L^1(\R).
\]
We call this case a ``relaxation layer".

 \noindent (b) If
$V>\int_0^\infty \xi g(\xi)\,d\xi$, we have
\[
\lim_{y\to\infty}F(y,\cdot)=M(-u_\infty,\cdot) \ \ \mathrm{in} \
L^1(\R).
\]
We call this case a ``shock layer".
\end{prop}
Let us first explain the terminology. In the first case of layer,
the first moment $\int \xi f\,d\xi$ of the given kinetic data at the
interface coincides with the first moment of the given Maxwellian at
infinity. In this case the process in the layer corresponds only to
a relaxation from  the kinetic data at the interface to its
Maxwellian. Note that in this case, the outgoing kinetic flow from
the layer is always 0.

In the second type of  layer, the value at $+\infty$ is a
Maxellian supported in $\{\xi\leq0\}$ with a flux stronger than the
in-going kinetic flux at the interface. This corresponds to a shock
front moving toward the interface. We can consider this kind of
layer as a shock profile which  ``sticks" to the interface. The
situations in the two cases of layer are very different. The first
set of initial conditions (\ref{initialrelax}) in  Theorem \ref{theolimitasymp}
ensures that, at the limit, the layer will be always of the
relaxation type. This implies in particular that, in this case, the
refined structure of the layer is not needed to get the boundary
condition. This is why we can still consider, in this case, weak limit
on the initial data. Those oscillations can destroy the kinetic
layer, but not the boundary condition. Things are completely
different for the second set of initial data (\ref{initialchoc}), which corresponds to
shock layer. Here, the precise structure of layer is needed to get
the correct boundary condition. This explains why this case is more
complicated. A precise control of the function $F_\eps$ is required.

Note that the uniqueness in Proposition \ref{theocouche} implies
that it does not exist any kinetic layer matching the Maxwellian
 $M(u_\infty,\xi)$ at $y=0$ with
$M(-u_\infty,\xi)$ at $y=+\infty$. This shows that we cannot have a
classical shock profile in a half space. This fact was already
proven by F.Golse in \cite{G}. This situation would correspond to
the critic situation between  {\it (a)} and {\it (b)}. Thus,
 we cannot go in a continuous way from a relaxation layer to a
shock layer.

The last remark has an interesting consequence. Consider the
following  initial data (non dependent on  $\eps$):
\begin{eqnarray*}
&&f^0(x,\xi)=M(u^+;\xi) \ \ \mathrm{for} \ \ x\leq0,\xi\in[-L,L],\\
&&f^0(x,\xi)=M(-u^+;\xi) \ \ \mathrm{for} \ \ x\geq0,\xi\in[-L,L].
\end{eqnarray*}
This initial condition verifies (\ref{initialrelax}), so Theorem
\ref{theolimitasymp} shows that $(\fe|_{\{x\leq0\}},u_\eps=\int
\fe\,d\xi|_{\{x\geq0\}})$ converges to the function (independent of
time):
\begin{eqnarray*}
&&f(t,x,\xi)=M(u^+;\xi) \ \ \mathrm{for} \ \
t\in\R^+,x\leq0,\xi\in[-L,L],\\
 &&u(t,x)=-u^+ \ \ \mathrm{for} \ \
t\in\R^+,x\geq0,\xi\in[-L,L].
\end{eqnarray*}
This corresponds to a steady shock localized at the interface $x=0$.
But as we have seen, there is no kinetic layer matching
$M(u^+;\cdot)$ to $-u^+$. Hence,  $F_\eps$ cannot converges to
(\ref{formel}). Physically, this means that at the $\eps$ scaling,
the shock is swept out from the interface. At the limit we only get
the relaxation layer.

In view of this pathology, the stability of the limit system could
seem surprising. The stability comes from the fact that, even if the
kinetic layer is not stable with respect to the boundary data, the
boundary condition $f(t,0-,\xi){\mathbf 1}_{\{\xi<0\}}$ itself is
stable.

\section{Preliminary results}

We begin with preliminary results on the model which will be useful
later. All those results are fairly standard and can be found, for example,
in \cite{Perthame}.  For every $u\in
[-L,L]$, we define $M(u,\cdot)\in L^\infty([-L,L])$ by:
\begin{equation}\label{defmax}
M(u,\xi)=\sign(u)\mathbf{1}_{\{0<\xi\sign u<|u|\}}.
\end{equation}
\begin{lemm}\label{calcul}
For every regular function $\phi\in C^\infty([-L,L])$ and
$u\in[-L,L]$ we have:
\[
\int_{-L}^{L}\phi'(\xi)M(u,\xi)\,d\xi=\phi(u)-\phi(0).
\]
\end{lemm}
 Let $N$ be an integer and  $\Ou$ be an open set of $\R^N$. Then for every
function $f\in L^\infty(\Ou;L^1([-L,L]))$ verifying $0\leq \sign
(\xi) f(\cdot,\xi)\leq 1$, we define $\M(f)\in
L^\infty(\Ou;L^1([-L,L]))$ by:
\begin{equation}\label{defmax2}
\M(f)(a,\xi)=M\left(\int_{-L}^{L}\!\!f(a,\zeta)\,d\zeta,\xi\right) \
\ a\in\Ou, \ \xi\in[-L,L].
\end{equation}
We have the following property:
\begin{lemm}\label{brenier} For every function $f\in
L^\infty([-L,L])$ compactly supported and such that
$0\leq\sign(\xi)f(\xi)\leq1$ and every convex function $\phi$ we
have:
\[
\int_{-L}^{L}\phi'(\xi)[\M f-f]\,d\xi\leq 0.
\]
Moreover, if for one strictly convex function $\phi$ we have the
equality:
\[
\int_{-L}^{L}\phi'(\xi)[\M f-f]\,d\xi= 0,
\]
then $f=\M f$.
\end{lemm}
\noindent {\bf Proof.} Notice that
\begin{eqnarray*}
\M f-f&\geq0&\mathrm{on} \ ]-L, \int f(\xi)\,d\xi[\\
&\leq0&\mathrm{on} \ ]\int f(\xi)\,d\xi,L[.
\end{eqnarray*}
Let us define
\[
h(\xi)=\int_{-L}^\xi(\M f-f)(\zeta)\,d\zeta.
\]
The function $h$ is null at $-L$, it is nondecreasing on $]-L, \int
f(\xi)\,d\xi[$ and non increasing on $]\int f(\xi)\,d\xi,L[$. From
Lemma \ref{calcul} and (\ref{defmax2}):
\[
\int_{-L}^{L}M\left(\int_{-L}^{L}f(\zeta)\,d\zeta;\xi\right)\,d\xi
=\int_{-L}^{L}f(\zeta)\,d\zeta.
\]
So $h$ is null at $L$. Hence $h(\xi)\geq0$ on $[-L,L]$. Finally,
integrating by parts give:
\[
\int_{-L}^{L}\phi'(\xi)(\M f-f)\,d\xi=-\int_{-L}^{L}
\phi"(\xi)h(\xi)\,d\xi\leq0
\]
if $\phi$ is convex. \\
Now consider one strictly convex function $\phi$. If the
inequality is an equality, then we have:
\[
\int_{-L}^{L} \phi"(\xi)h(\xi)\,d\xi=0,
\]
which implies that $h$ is null for almost every $\xi\in[-L,L]$.
Because of the definition of $h$, this implies that $\M f=f$ for
almost every $\xi\in[-L,L]$. \qquad \qed
\begin{lemm}\label{simple}
Let $g_n\in L^\infty(\Ou;L^1([-L,L]))$, $0\leq \sign(\xi)g_n\leq1$
and $u\in L^1(\Ou)$ be such that $g_n$ converges weakly to
$M(u,\xi)$ in $L^\infty w*$. Then, the convergence holds strongly in
$L^1_{\mathrm{loc}}(\Ou\times[-L,L])$.
\end{lemm}
\noindent {\bf Proof.} Let $\Omega$ be a compact set of
$\Ou\times[-L,L]$. $g_n$ is bounded in $L^\infty$ so in
$L^2(\Omega)$. The convergence holds weakly in $L^2w$. But:
\begin{eqnarray*}
&&\int_\Omega
|M(u(x),\xi)|^2\,dx\,d\xi\leq\underline{\lim}\int_\Omega
|g_n|^2\,dx\,d\xi\\
&&\leq\overline {\lim}\int_\Omega
|g_n|^2\,dx\,d\xi\\
&&\leq \lim\int_\Omega \sign(\xi)g_n(x,\xi)\,dx\,d\xi
=\int_\Omega\sign(\xi)M(u(x),\xi)\,dx\,d\xi.
\end{eqnarray*}
Since $|M(u,\xi)|^2=\sign(\xi)M(u,\xi)$, we conclude that the
$L^2(\Omega)$ norm of $g_n$ converges to the $L^2(\Omega)$ norm of
$M(u,\xi)$. Hence $g_n$ converges strongly in $L^2(\Omega)$ and so
in $L^1_{\mathrm{loc}}(\Ou\times[-L,L])$.\qed
\begin{lemm}\label{cadecroit}
Consider two functions  $f,g\in L^1([-L,L])$ verifying the
compatibility conditions $0\leq \sign(\xi)f(\xi)\leq 1$ and $0\leq
\sign(\xi)g(\xi)\leq 1$. Then:
\begin{equation}\label{inequal}
\int_{-L}^L \mathbf{1}_{\{f\geq g\}} \left[(\M f(\xi)-\M
g(\xi))-(f(\xi)-g(\xi)) \right]\,d\xi\leq0.
\end{equation}
Moreover, if  this quantity is equal to 0 then  $\sign(f-g)$ is
constant  on $[-L,L]$. Especially, in this case, $\sign(\M f-\M
g)=\sign(f-g)$.
\end{lemm}
By $\sign(f-g)$ is
constant  on $[-L,L]$, we mean that either $f(\xi)\leq g(\xi)$ for all $\xi\in (-L,L)$, or $f(\xi)\geq g(\xi)$ for all $\xi\in (-L,L)$.

\noindent {\bf Proof.} If $\int f\,d\xi\leq\int g\,d\xi$ then $\M
f-\M g\leq0$ for every $\xi\in[-L,L]$ and the inequality is true.
If, moreover,  (\ref{inequal}) is an equality, then %considering the
%$\xi$ such that $\M f-\M g\neq0$ and $\M f-\M g=0$,
$$
\int_{-L}^L\mathbf{1}_{\{f\geq
g\}}[f(\xi)-g(\xi)]\,d\xi=\int_{-L}^L\mathbf{1}_{\{f\geq g\}}[\M
f(\xi)-\M g(\xi)]\,d\xi\leq 0,
$$
and so $f\leq g$ on $[-L,L]$.

Assume now that  $\int f\,d\xi\geq\int g\,d\xi$. Then $\M f-\M
g\geq0$ for every $\xi\in[-L,L]$ and:
\begin{eqnarray*}
&&\int_{-L}^L \mathbf{1}_{\{f\geq g\}}(\M f(\xi)-\M
g(\xi))\,d\xi\leq
\int_{-L}^L(\M f(\xi)-\M g(\xi))\,d\xi\\
&&\qquad\qquad=\int_{-L}^L(f(\xi)-g(\xi))\,d\xi\\
&\leq&\int_{-L}^L\mathbf{1}_{\{f\geq g\}}(f(\xi)-g(\xi))\,d\xi+
\int_{-L}^L\mathbf{1}_{\{f\leq g\}}(f(\xi)-g(\xi))\,d\xi\\
&\leq&\int_{-L}^L\mathbf{1}_{\{f\geq g\}}(f(\xi)-g(\xi))\,d\xi,
\end{eqnarray*}
Hence the inequality is verified. If (\ref{inequal}) is an equality,
then
\[
\int_{-L}^L(f(\xi)-g(\xi))\,d\xi=\int_{-L}^L\mathbf{1}_{\{f\geq
g\}}(f(\xi)-g(\xi))\,d\xi,
\]
so $f\geq g$ on $[-L,L]$, which ends the proof.\qquad\qed

\section{Study of the layer}

This section is devoted to the proof of Proposition
\ref{theocouche}.
 Notice that it is necessary for $(V,g)$ to
be in $\C$ in order to have a solution to (\ref{Milne}). Indeed,
if $F$ is solution to (\ref{Milne}), then integrating $\xi
F(0,\xi)$ with respect to $\xi$ we find:
\[
V=\int_{-L}^{L} \xi F(0,\xi)\,d\xi=\int_0^L \xi
g(\xi)\,d\xi+\int_{-L}^0\xi F(0,\xi)\,d\xi\geq\int_0^L \xi
g(\xi)\,d\xi.
\]
Indeed $\int_{-L}^0\xi F(0,\xi)\,d\xi\geq0$ because of the condition
$0\leq \sign{\xi}F(y,\xi)$. \vskip3mm

As said in the introduction, in the case of a relaxation layer, the
value $F(y=0,\xi\leq0)$ used in the limit coupled problem is 0.
Hence, in this case, the precise structure of the layer is not
needed. In some cases the layer function $\Fe$ does not converge to
the corresponding layer solution. But this can occur only in the
case of the relaxation layer. And we will show that, fortunately, we
still recover the limit problem since the coupling value
$f(t,0-,\xi\leq0)=0$ is verified, which is the same that the value
provided by the layer solution.

 In the contrary, the
precise structure of the shock layer is needed to defined the limit
problem. Fortunately, in this case the layer function $\Fe$
converges to the layer solution. \vskip0.3cm

We decompose the proof into three lemmas. The first one is about the
properties of solutions to the layer problem, the second one about
the uniqueness and the third one about the existence of the
solution.

\begin{lemm}\label{prop_proprieteslayer}
Let $(V,g)\in\C$. We set  $u_\infty\in [0,L[$ such that:
\[
V=\frac{u_\infty^2}{2}.
\]
Assume that $F\in L^\infty(\R^+\times]-L,L[)$ verifies
$0\leq\sign(\xi)F\leq1$ and is solution to the layer problem
(\ref{Milne}). Then, we have the following property depending on the
boundary data.

\noindent (a) Relaxation layer: \\
if $V=\int_0^L \xi g(\xi)\,d\xi$ then:
\[
\lim_{y\to\infty}F(y,\cdot)=M(u_\infty,\xi) \ \ \mathrm{in} \
L^1(]-L,L[).
\]
Moreover $F(y,\xi)=0$ for every $y>0$ and
$\xi\in]-L,0[$.\vskip0.3cm

\noindent (b) shock layer:\\
if $V>\int_0^L \xi g(\xi)\,d\xi$ then:
\[
\lim_{y\to\infty}F(y,\cdot)=M(-u_\infty,\xi) \ \ \mathrm{in} \
L^1(]-L,L[).
\]
\end{lemm}
\noindent{\bf Proof.} We divide the proof into several
steps.\vskip2mm

{\it (i) First Property of relaxation layers.} We consider the case
of relaxation layers, namely $V=\int_0^L \xi g(\xi)\,d\xi$. We have
\[
V=\int_0^L \xi g(\xi)\,d\xi+\int_{-L}^0 \xi F(0,\xi)\,d\xi.
\]
But $\sign(\xi) F\geq0$ Hence $F(0,\xi)=0$ for $\xi<0$. Since $\M
F$ is non positive for $\xi\leq0$, equation (\ref{Milne}) gives
for $\xi<0$:
\[
\dy (\xi F)=\M F-F\leq -\frac{\xi F}{\xi}.
\]
So, for every $y>0,\xi<0$:
\[
\xi F(y,\xi)\leq \xi F(0,\xi) e^{-y/\xi}=0.
\]
But $\sign(\xi F)\geq0$, hence:
\begin{equation}
F(y,\xi)=0 \ \ \mathrm{for} \  y\in\R^+, \ \xi\in]-L,0[.
\end{equation}
{\it (ii) Preliminary result for shock layers.} We consider the case
of shock layer, namely $V>\int_0^L \xi g(\xi)\,d\xi$. Multiplying
equation (\ref{Milne}) by $\mathbf{1}_{\{0\leq\xi\}}$ we find
(thanks to Lemma \ref{brenier}):
\[
\int_0^L \xi F(y,\xi)\,d\xi
\]
is non increasing. But
\[
V=\int_{-L}^0 \xi F(y,\xi)\,d\xi+\int_0^L \xi F(y,\xi)\,d\xi.
\]
Hence for $y>0$:
\begin{equation}
\int_{-L}^0\xi F(y,\xi)\,d\xi\geq V-\int_0^L \xi g(\xi)\,d\xi>0.
\end{equation}
\vskip2mm \noindent{\it (iii) Vanishing entropy up to a
subsequence.}  Multiplying the first equation of (\ref{Milne}) by $(-\xi)$ and
 integrating it on $[0,y]\times(-L,L)$, we find
\begin{eqnarray*}
&&\qquad\int_0^y\int_{-L}^L (-\xi)[\M
 F(y,\xi)-F(y,\xi)]\,d\xi\,dy\\
&&\leq \int_{-L}^L\xi^2 F(0,\xi)\,d\xi-\int_{-L}^L\xi^2 F(y,\xi)\,d\xi\\
&&\leq 2\frac{L^3}{3},
\end{eqnarray*}
since $|F|\leq 1$.
Thanks to Lemma \ref{brenier}, $\int (-\xi)[\M
 F(y,\xi)-F(y,\xi)]\,d\xi$ is nonnegative, so it is bounded in $L^1(\R^+)$ as a function of $y$.
 Hence,
there exists a sequence
 $y_n\to\infty$ such that $\int (-\xi)[\M
 F(y_n,\xi)-F(y_n,\xi)]\,d\xi$  converges to 0.

\vskip2mm \noindent{\it (iv) Convergence of $\M F(y_n,\cdot)$ up to
a subsequence.} Let us denote $u(y)=\int_{-L}^{L}F(y,\xi)\,d\xi$.
Thanks to the definition to $\M$:
\[
\frac{u^2(y_n)}{2}=\int_{-L}^{L}\xi \M F(y_n,\xi)\,d\xi.
\]
Using the definition of $V$ and the result {\it (iii)}, we find that
\[
\frac{u^2(y_n)}{2}-V=\int_{-L}^{L}\xi [\M
F(y_n,\xi)-F(y_n,\xi)]\,d\xi
\]
 converges to 0 when $n$ goes to infinity. Hence, up to a subsequence, $u(y_n)$
converges to $u_\infty$ or to $-u_\infty$. This means that, up to a
subsequence,  $\M F(y_n,\cdot)$ converges in $L^1(]-L,L[)$ to
$M(u_\infty,\cdot)$ or to $M(-u_\infty,\cdot)$, when $n$ goes to
infinity. \vskip2mm\noindent{\it (v) Convergence of $F(y_n,\cdot)$.}
 Since $\|F\|_{L^\infty}\leq1$, extracting a subsequence from the
 above one if necessary, $F(y_n,\cdot)$ converges weakly in $L^\infty*$ to a function
$F_\infty$ such that $0\leq \sign(\xi)F_\infty\leq1$.  Thanks to
{\it (iii)} and the strong convergence of $\M F(y_n,\cdot)$ stated
in {\it (iv)} we find at the limit:
\[
\int_{-L}^L\xi[ F_\infty(\xi)-\M F_\infty(\xi)]\,d\xi=0.
\]
Thanks to Lemma \ref{brenier}, this implies that $F_\infty=\M
F_\infty$. Hence $F_\infty(\xi)$ is $M(u_\infty,\cdot)$ or
$M(-u_\infty,\cdot)$. Thanks to Lemma \ref{simple}, this implies
that  the convergence holds strongly in $L^1(]-L,L[)$. But thanks to
{\it (i)}, the limit cannot be  $M(-u_\infty,\cdot)$ (at least if
$u_\infty\neq0$) in the case of relaxation layers. So in this case
the limit is $M(u_\infty,\cdot)$. In the case of shock layers,
thanks to {\it (ii)}, $M(u_\infty,\cdot)$ cannot be the limit. Hence
in this case, the limit is  $M(-u_\infty,\cdot)$. By the uniqueness
of the limit,  the entire sequence $F(y_n,\xi)$ converges strongly
in $L^1(\R)$ to $M(u_\infty,\cdot)$ in the case of relaxation layer
and to $M(-u_\infty,\cdot)$ in the case of shock layer.
\vskip2mm\noindent{\it (vi) Convergence for $y\to\infty$.} Consider
a monotonic function $h$. Multiplying the first equation to
(\ref{Milne}) by this function and integrating with respect to
$\xi$, we find thanks to Lemma \ref{brenier} that $\int \xi
h(\xi)F(y,\xi)\,d\xi$ is monotonic  with respect to $y$. Notice that
since $F$ is bounded by 1 and compactly supported in $[-L,L]$, this
function is bounded and is converging to a constant when $y$ goes to
infinity. Thanks to {\it (v)} the limit is $\int \xi
h(\xi)M(u_\infty,\xi)\,d\xi$ in the case of relaxation layer and
$\int \xi h(\xi)M(-u_\infty,\xi)\,d\xi$ in the case of shock layer.
Since every regular function is the difference between two
nondecreasing functions, this convergence holds true for every
regular function $h$. This implies that for every $\eta>0$, The
whole family $F(y,\cdot)$ converges in $\D'(]\eta,L[)$ and in
$\D'(]-L,-\eta[)$ to the corresponding limit function. But since
those functions are uniformly bounded in $L^\infty(]-L,L[)$, finally
the convergence holds in $\D'(]-L,L[)$. But the limit function is an
equilibrium function, hence, thanks to Lemma \ref{simple}, the
convergence holds strongly in $L^1(]-L,L[)$.\qed
 \vskip3mm
 Let us now consider the uniqueness of solution to the
layer problem.
\begin{lemm}\label{lemm_uniqueness}
For every $(V,g)\in\C$ there exists at most one  solution $F$ to (\ref{Milne}) verifying
  $F\in
L^\infty(\R^+\times]-L,L[)$ with $0\leq\sign(\xi)F\leq1$.
\end{lemm}
\noindent{\bf Proof.} Consider $F_1,F_2$ two solutions to problem
(\ref{Milne}) for the same condition values $(V,g)\in\C$.
Multiplying the difference of the first equations of (\ref{Milne})
for $F_1$ and $F_2$ by $\sign(F_1-F_2)$, and integrating in $\xi$, we find:
\[
\dy\int_{-L}^{L}\xi|F_1-F_2|\,d\xi=\int_{-L}^{L}\sign(F_1-F_2)[\M
F_1-\M F_2-F_1+F_2 ]\,d\xi.
\]
Thanks to Lemma \ref{cadecroit} this quantity is non-positive. Thanks
to Lemma \ref{prop_proprieteslayer}, we have:
\[
\left. \int_0^L\xi|F_1-F_2|\,d\xi\right|_{y=0}=0,
\]
and:
\[
\lim_{y\to \infty }\int_{-L}^{L}\xi|F_1-F_2|\,d\xi=0.
\]
Hence:
\[
\int_0^\infty\int_{-L}^{L}\sign(F_1-F_2)[\M F_1-\M F_2-F_1+F_2
]\,d\xi\,dy=0
\]
and:
\[
\left.\int_{-L}^{0}\xi|F_1-F_2|\,d\xi\right|_{y=0}=0.
\]
Thanks to Lemma \ref{cadecroit} this implies that
$\partial_\xi(\sign(F_1-F_2))=0$. Let us  fix a $\xi<0$. We denote
$\Omega_\xi=\{y|F_1(y,\xi)>F_2(y,\xi)\}$. Notice that since $|\dy
(F_1-F_2)|\leq4/|\xi|$, $(F_1-F_2)(\cdot,\xi)$ is continuous and so
$\Omega_\xi$ is an open subset.  Assume that it is not empty. Denote
$y>0$ one of its elements and $y_0=\inf\{z<y|
]z,y[\subset\Omega_\xi\}$. For every $z\in]y_0,y[$ and every
$\zeta\in]-L,L[$ we have: $F_1(z,\xi)> F_2(z,\xi)$, hence we have
$F_1(z,\zeta)\geq F_2(z,\zeta)$ too. This leads to:
\[
\xi \dy (F_1-F_2)+(F_1-F_2)\geq \M F_1-\M F_2\geq0, \qquad
\mathrm{on \ } ]y_0,y[.
\]
Hence, we find that $y_0=0$ and $(F_1-F_2)(0,\xi)\geq
(F_1-F_2)(y,\xi)e^{y/\xi}>0 $. This gives a contradiction. Hence
$\Omega_\xi$ is empty and so $F_1\leq F_2$ for $\xi<0$. Exchanging
the indices gives:
\[
F_1(y,\xi)=F_2(y,\xi) \qquad y\in\R^+, \ \xi\in[-L,0].
\]
Hence, for every $y>0$
\[
\int_0^L \xi|F_1-F_2|(y,\xi)\,d\xi=\int_{-L}^L
\xi|F_1-F_2|(y,\xi)\,d\xi\leq \int_{-L}^L
\xi|F_1-F_2|(0,\xi)\,d\xi=0.
\]
This implies that $(F_1-F_2)(y,\xi)=0$ for positive $\xi$ too. This
ends the proof.\qed
 \vskip3mm
 Let us now show the existence of the solution to
the layer problem which ends the proof of Proposition
\ref{theocouche}.
\begin{lemm}\label{prop_existencelayer}
For every condition values $(V,g)\in\C$, there exists a  solution
$F\in L^\infty(\R^+\times]-L,L[)$ with $0\leq\sign(\xi)F\leq1$,
solution to the layer problem (\ref{Milne}).

\end{lemm}
\noindent{\bf Proof.} We divide the proof into several parts.
 \vskip2mm\noindent{\it (i) Construction of approximated solutions.}\\
In order to show the existence of the solution to layer problem
(\ref{Milne}), we use the method of Golse \cite{G}. Let
$(V,g)\in\C$. We construct by induction $F_n$ in the following way.
We set $F_1=0$ for the case (a) and $F_1=M(-u_\infty,\xi)$ for the
case (b) and for $n>1$:
\begin{eqnarray*}
F_{n+1}(y,\xi)&=&g(\xi)e^{-\frac{y}{\xi}}+\int_0^y\frac{1}{\xi}\M
F_n(z,\xi)e^{\frac{z-y}{\xi}}\,dz \ \ \mathrm{for}
 \ y\geq0,\xi\in]0,L[,\\
&=&-\int_y^\infty\frac{1}{\xi}\M F_n(z,\xi)e^{\frac{z-y}{\xi}}\,dz \
\ \mathrm{for} \ y\geq0,\xi\in]-L,0[.
\end{eqnarray*}
We can check  that $F_n$ verifies $0\leq \sign(\xi)F_n\leq1$:
\begin{equation}\label{Milne_approchee}
\left\{
\begin{array}{l}
\dis{\xi\dy F_{n+1}=\M F_n-F_{n+1}\qquad y>0, \ \xi\in]-L,L[,}\\[3mm]
\dis{F_{n+1}(y=0,\xi\geq0)=g\qquad \xi\in]0,L[.}%\\[3mm]
%\dis{F_{n+1}(y,\xi\leq0)=0,}
\end{array}
\right.
\end{equation}

\noindent{\it (ii) Convergence when $n\to\infty$.}\vskip2mm Let us
show by induction that for every fixed point
$(y,\xi)\in\R^+\times]-L,L[$, the sequence $\{F_n(y,\xi)\}$ is non
decreasing. First check that
\[
F_2(y,\xi)-F_1(y,\xi)=g(\xi)e^{-\frac{y}{\xi}}\mathbf{1}_{\{\xi\geq0\}}\geq0.
\]
Assume now that $F_n(y,\xi)\geq F_{n-1}(y,\xi)$ for every
$(y,\xi)\in\R^+\times]-L,L[$. Then
\begin{eqnarray*}
&&F_{n+1}(y,\xi)-F_n(y,\xi)\\
&&=\int_0^y\frac{1}{\xi}\left[\M F_n(z,\xi)-\M
F_{n-1}(z,\xi)\right]e^{\frac{z-y}{\xi}}\,dz
\geq0 \ \ \mathrm{for} \ \xi\in]0,L[,\  y\geq0,\\
&=&-\int_y^\infty\frac{1}{\xi}\left[\M F_n(z,\xi)-\M
F_{n-1}(z,\xi)\right]e^{\frac{z-y}{\xi}}\,dz\geq0 \ \ \mathrm{for} \
\xi\in]-L,0[,\ y\geq0.
\end{eqnarray*}
By this procedure we have shown that the sequence $\{F_n(y,\xi)\}$
is non decreasing. Moreover it is bounded by 1, hence it converges
 almost everywhere to a function $F(y,\xi)$ with
 $0\leq\sign(\xi)F(y,\xi)\leq1$. By Lebesgue's Theorem, $u_n(y)=\int_{-L}^L
 F_n(y,\xi)\,d\xi$ converges almost everywhere to $u_F(y)=\int_{-L}^L
 F(y,\xi)\,d\xi$. So thanks to the definition to $\M$,
 $\M F_n$ converges strongly to $\M F$ in
 $\Lloc^1(\R^+\times]-L,L[)$. Passing to the limit in the first
 equation to (\ref{Milne_approchee}) shows that $F$ is solution to the first equation to
 (\ref{Milne}). The first equation of (\ref{Milne_approchee})
 shows that $\dy(\xi F_n)$ is bounded in
 $L^\infty(\R^+\times]-L,L[)$. Since $F_n\in
 L^\infty(\R^+\times]-L,L[)$ and $L^\infty(]-L,L[)$ is compactly  imbedded in
 $H^{-1}(]-L,L[)$, from Aubin's Theorem we find that the
 convergence holds in $C^0(\R^+;H^{-1}(]-L,L[))$. Hence we
 retrieve at the limit $n\to\infty$
 \[
F(0,\xi\geq0)=g\qquad\qquad y>0,\ \xi\in]0,L[.
 \]
 \noindent{\it (iii) Flux condition inside the layer.}\vskip2mm
 We have now to show the last condition of (\ref{Milne}).
 Let us first consider the relaxation case:\\
 (a) $V=\int_0^L \xi g(\xi)\,d\xi$. First, integrating the
 first equation of (\ref{Milne}) with respect to $\xi$
we find that
\[
\int_{-L}^L\xi F(y,\xi)\,d\xi=\int_{-L}^L\xi F(0,\xi)\,d\xi\qquad
y\geq0.
\]
Since $F_1(y,\xi)=0$ for every $(y,\xi)\in\R^+\times]-L,0[$ and
$F_n$ is non decreasing, we have that $F(y,\xi)\geq 0$ for
$(y,\xi)\in\R^+\times]-L,0[$. But at the limit $0\leq
\sign(\xi)F(y,\xi)\leq1$ hence $F(y,\xi)=0$ for
$(y,\xi)\in\R^+\times]-L,0[$. So
\[
V=\int_0^L \xi g(\xi)\,d\xi=\int_{-L}^L\xi
F(0,\xi)\,d\xi=\int_{-L}^L\xi F(y,\xi)\,d\xi\qquad y>0.
\]
Consider now the shock case:\\
(b) $V>\int_0^L \xi g(\xi)\,d\xi$. By induction we show that for
every $n>0$, $F_n(y,\cdot)$ converges to $M(-u_\infty,\cdot)$ in
$L^1(]-L,L[)$ when $y\to\infty$. This is obviously true for $n=1$
and if it is true for $n$ then, $\M F_n(y,\cdot)$ converges to
$M(u_\infty,\cdot)$ in $L^1(\R)$ too. The results follows from the
definition to $F_{n+1}$. Integrating the first equation of
(\ref{Milne_approchee}) with respect to $\xi$ and using the
nondecreasing property of $\{F_n(y,\xi)\}$ we find:
\begin{eqnarray*}
\dy\int_{-L}^L \xi F_n(y,\xi)\,d\xi&=&\int_{-L}^L (\M
F_{n-1}(y,\xi)-F_n(y,\xi))\,d\xi\\
&\leq&\int_{-L}^L (\M F_{n}(y,\xi)-F_n(y,\xi))\,d\xi=0.
\end{eqnarray*}
Hence using the limit at $+\infty$, we find for every $y$:
\[
\int_{-L}^L\xi F_n(y,\xi)\,d\xi\geq \int_{-L}^L\xi
M(-u_\infty,\xi)\,d\xi=V.
\]
Passing to the limit we find:
\[
\int_{-L}^L\xi F(y,\xi)\,d\xi=V_0\geq V.
\]
In particular $V_0>\int\xi g(\xi)\,d\xi$ hence $F$ is a layer of
type (b). Thanks to Lemma \ref{prop_proprieteslayer}:
\[
\lim_{y\to\infty}F(y,\cdot)=M(-u_0,\cdot),
\]
where $u_0\in[0,L[$ verifies $V_0=u_0^2/2$. Since $F\geq F_1$, we
have $F(y,\xi)=0$ for $\xi\leq -u_\infty$. Hence $u_0\leq u_\infty$.
So finally $V_0\leq V$ and so $V_0=V$.
 \qquad\qed\vskip3mm
 We finish this section with the following Proposition which gives
a property of good confinement of the layers.
\begin{prop}\label{prop_confinement}
For any  $(V,g)\in\C$, the function $F$ solution to  (\ref{Milne})
verifies
$$
\int_0^\infty\int_{-L}^L
|F(y,\xi)-F_\infty(\xi)|\,d\xi\,dy<\infty,
$$
where
$$
F_\infty=\lim_{y\to\infty}F(y,\cdot).
$$
\end{prop}
\noindent{\bf Proof.}
If $V=0$, then $F$ is identically 0 and the result holds.
If not, we have $V=(u_\infty)^2/2>0$. From Lemma \ref{prop_proprieteslayer},
$$
F_\infty(\xi)=M(u_*,\xi),
$$
with $u_*=\pm u_\infty$. Let $u(y)=\int_{-L}^LF(y,\xi)\,d\xi$. We have
$$
\lim_{y\to\infty} u(y)=u_*,
$$
and, thanks to Lemma \ref{brenier},
$$
(u(y))^2-(u_*)^2=(u(y))^2-2V=2\int_{-L}^L\xi(\M F(y,\xi)-F(y,\xi))\,d\xi\leq 0.
$$
Multiplying the layer equation by $2\xi$ and integrating in $\xi$, we find, for $0<y<\infty$,
$$
2\partial_y\int_{-L}^L \xi^2F(y,\xi)\,d\xi=2\int_{-L}^L \xi(\M F(y,\xi)-F(y,\xi))\,d\xi=(u(y))^2-(u_*)^2.
$$
Integrating in $y$ gives that
$$
\int_0^\infty [(u(y))^2-(u_*)^2]\,dy<\infty.
$$
Note that $|u_*-u(y)|=((u_*)^2-(u(y))^2)/|u_*+u(y)|$, and for $y$ large enough we have $|u_*+u(y)|>|u_*|>0$. Hence
$$
\int_0^\infty |u(y)-u_*|\,dy<\infty.
$$
Using $\M F(y,\xi)=M(u(y),\xi)$, we find
$$
\int_0^\infty\int_{-L}^L|\M F(y,\xi)-F_\infty(\xi)|\,d\xi\,dy=\int_0^\infty |u(y)-u_*|\,dy<\infty.
$$
Finally, $|F(y,\xi)-F_\infty(\xi)|=\sign(\xi-u_*)(F(y,\xi)-F_\infty(\xi))$ and
\begin{eqnarray*}
&&\qquad \sign(\xi-u_*)\xi\partial_y F(y,\xi)\\
&&=\sign(\xi-u_*)(\M F(y,\xi)-F_\infty(\xi))-\sign(\xi-u_*)(F(y,\xi)-F_\infty(\xi))\\
&&=\sign(\xi-u_*)(\M F(y,\xi)-F_\infty(\xi))-|F(y,\xi)-F_\infty(\xi)|.
\end{eqnarray*}
Integrating this expression in $y$ and $\xi$ gives
$$
\int_0^\infty\int_{-L}^L|F(y,\xi)-F_\infty(\xi)|\,d\xi\,dy\leq \int_0^\infty\int_{-L}^L|\M F(y,\xi)-F_\infty(\xi)|\,d\xi\,dy+2L^2<\infty.
$$
This ends the proof.\qquad \qed

\section{Well-posedness  of the limit problem}

This section is devoted to the proof of Theorem \ref{theoprolimite}
which shows that the coupled system
(\ref{limitefluide}) (\ref{limitecinetique})(\ref{limitecouche}) is
well-posed and stable with respect to the initial values.

A careful reader could be surprised that the limit problem is stable
with respect to the initial conditions. Indeed,  in the layer,
$F_1-F_2$ cannot be controlled by the initial conditions.  However,
notice that the coupling in the limit problem depends only on the
values at $y=0$ of the layer. And this value is continuous with
respect  to the conditions of the problem.

In order to obtain the existence of a solution to the limit problem,
we use the fixed point Theorem of Schauder applied on the trace of
$u$ at the interface $x=0$. The result relies on the stability  of
some quantities  at $x=0$.

We first show a stability result in the kinetic domain.
\begin{lemm}\label{lemm_kineticaubord}
Let  $f_1^0$, $f_2^0\in L^\infty((-\infty,0)\times]-L,L[)$ and $g_1$
$g_2\in L^\infty((0,\infty)\times]-L,0[)$. %with values between $-1$ and 0.
Consider $f_1$ $f_2\in
L^\infty((0,\infty)\times(-\infty,0)\times]-L,L[)$ solutions  for $i=1,2$
to
\begin{equation*}
\left\{
\begin{array}{l}
\dt f_i+\xi\dx f_i=\M f_i-f_i, \ \ \mathrm{for} \
t\in\R^+,x\in\R^-,\xi\in]-L,L[,\\[3mm]
f_i(t=0,x,\xi)=f^0_i(x,\xi)  \ \ \mathrm{for} \ x\in\R^-,\xi\in]-L,L[,\\[3mm]
f_i(t,x=0-,\xi\leq0)=g_i(t,\xi\leq0) \ \ \mathrm{for} \
t\in\R^+,\xi\in]-L,0[.
\end{array}
\right.
\end{equation*}
Then, we have, for every $t>0$
\begin{eqnarray*}
&&\int_{-\infty}^0\int_{-L}^L|f_1(t)-f_2(t)|\,d\xi\,dx+\int_0^t\int_{0}^L\xi|f_1(s,0-,\xi)-f_2(s,0-,\xi)|\,d\xi\,ds\\
&&\qquad\leq
\int_{-\infty}^0\int_{-L}^L|f_1^0-f_2^0|\,d\xi\,dx+\int_0^t\int_{-L}^0(-\xi)|g_1(s,\xi)-g_2(s,\xi)|\,d\xi\,ds.
\end{eqnarray*}
Moreover, if $f_\eps^0$ converges strongly in $L^1_\mathrm{loc}$ to
$f^0$ and $g_\eps$ converges weakly to $g$, then
the function $f_\eps(t,0-,\xi\geq0)$ converges STRONGLY to $f(t,0-,\xi\geq0)$ in
$L^1_\mathrm{loc}((0,\infty)\times]0,L[)$, where $f$ is the solution
with initial data $f^0$ and boundary data $g$.
\end{lemm}
\noindent{\bf Proof.} Thanks to Lemma \ref{cadecroit}, we have
$$
\dt \int_{-L}^L|f_1-f_2|\,d\xi+\dx \int_{-L}^L\xi|f_1-f_2|\,d\xi\leq 0\ \ \mathrm{for} \
t\in\R^+,x\in\R^-.
$$
Integrating in $(t,x)$ gives the first result. For the second
one, first notice that, thanks to averaging lemmas (introduced first in \cite{lm1}. See \cite{lm2} for our case.), $\int
f_\eps(t,x,\xi)\,d\xi$ converges strongly in $L^1_{\mathrm{loc}}((0,\infty)\times\R^-)$. Hence $\M f_\eps$
converges strongly to $\M f$ in $L^1_{\mathrm{loc}}((0,\infty)\times\R^-\times(-L,L))$. Then the Duhamel formula gives for
$\xi>0$:
$$
f_\eps(t,0-,\xi)=f^0_\eps(-t\xi,\xi)e^{-t}+\int_0^t e^{s-t}\M
f_\eps(s,(s-t)\xi,\xi)\,ds,
$$
which gives the result.\qed

We now give a similar result for the coupling of the fluid part with the layer.
\begin{lemm}\label{lemm_coucheetbord}
For any $g\in L^\infty(\R^+\times(0,L))$ with
$0\leq g(t,\xi)\leq1$, and $u^0\in L^\infty(0,\infty)$, $-L\leq u^0(x)\leq L$, there exists a unique solution $(u,F)$ to
\begin{equation*}
\left\{
\begin{array}{l}
\dis{\dt u+\dx\frac{u^2}{2}=0,} \ \ \mathrm{for} \
t\in\R^+,x\in\R^+, \\[3mm]
u(t=0,x)=u^0(x), \ \
\mathrm{for} \ x\in\R^+,
\\[3mm]
\dis{u(t,x=0+)\stackrel{BLN}{:=}\sqrt{2\int_0^\infty\xi
g(t,\xi)\,d\xi}, \ \ \mathrm{for} \ \ t\in\R^+,}
\end{array}
\right.
\end{equation*}

\begin{equation*}
\left\{
\begin{array}{l}
\xi\dy F=\M F-F, \ \ \mathrm{for} \
t\in\R^+,y\in\R^+,\xi\in(-L,L),\\[3mm]
F(t,y=0+,\xi\geq0)=g(t,\xi), \ \ \mathrm{for} \
t\in\R^+,\xi\in(0,L),\\[3mm]
\dis{\int_{-\infty}^{+\infty}\xi
F(t,y,\xi)\,d\xi=\frac{u(t,0+)^2}{2}} \ \ \mathrm{for} \
t\in\R^+,y\in\R^+.
\end{array}
\right.
\end{equation*}
Moreover, if $(u_1,F_1)$ and $(u_2,F_2)$ are two such solutions associated to
$g_1$, $u_1^0$ (respectively $g_2$, $u_2^0$), then, for any $t>0$,
\begin{eqnarray*}
&&\int_0^\infty|u_1(t,x)-u_2(t,x)|\,dx+\int_0^t\int_{-L}^0|\xi||F_1(s,0+,\xi)-F_2(s,0+,\xi)|\,d\xi\,ds\\
&&\qquad\leq\int_0^\infty|u^0_1(x)-u^0_2(x)|\,dx+\int_{0}^t\int_0^L\xi|g_1(s,\xi)-g_2(s,\xi)|\,d\xi\,ds.
\end{eqnarray*}
\end{lemm}
\noindent{\bf Proof.}
For any $g$, we can construct a solution to the initial-boundary problem of the Burgers equation (see \cite{otto} or \cite{Vovelle}, for instance). This solution constructed, the strong trace theorem \cite{V1} gives a meaning to $u(t,0+)$. We can then consider the solution of the layer constructed in Proposition \ref{theocouche}.

Consider now two solution $(u_1,F_1)$ and $(u_2,F_2)$ of this problem.
From the Kruzkov's theory of  the Burgers equation (\cite{Kruzkov}), and using the strong trace theorem,
we find for any $t>0$
\begin{eqnarray*}
&&\qquad\qquad\dt \int_0^\infty|u_1(t,x)-u_2(t,x)|\,dx\\
&&\leq \sign(u_1(t,0+)-u_2(t,0+))(u_1(t,0+)^2-u_2(t,0+)^2)/2.
\end{eqnarray*}
Using Lemma \ref{cadecroit}, we find that
$$
\dy\int_{-L}^L\xi|F_1(t,y,\xi)-F_2(t,y,\xi)|\,d\xi\leq0.
$$
Integrating in $y$ between 0 and $\infty$ and using
Proposition \ref{theocouche} give
\begin{eqnarray*}
&&\qquad\qquad\int_{-L}^0|\xi||F_1(t,0,\xi)-F_2(t,0,\xi)|\,d\xi\\
&&\leq \int_0^L\xi|g_1(t,\xi)-g_2(t,\xi)|\,d\xi-\int_{-L}^L\xi|M(u^\infty_1(t,\xi),\xi)-M(u^\infty_2(t,\xi),\xi)|\,d\xi,
\end{eqnarray*}
where
$$
\lim_{y\to\infty}F_1(t,y,\xi)=M(u^\infty_1(t),\xi),\qquad \lim_{y\to\infty}F_2(t,y,\xi)=M(u^\infty_2(t),\xi).
$$
Note that
\begin{eqnarray*}
&&\qquad\int_{-L}^L\xi|M(u^\infty_1(t,\xi),\xi)-M(u^\infty_2(t,\xi),\xi)|\,d\xi\\
&&=\sign(u^\infty_1(t)-u^\infty_2(t))(u^\infty_1(t)^2-u^\infty_2(t)^2)/2,
\end{eqnarray*}
and
$$
u_1^2={u^\infty_1}^2,\qquad u_2^2={u^\infty_2}^2.
$$
Putting those results together gives
\begin{eqnarray}
&&\nonumber\dt\int_0^\infty|u_1(t,x)-u_2(t,x)|\,dx+\int_{-L}^0|\xi||F_1(t,0+,\xi)-F_2(t,0+,\xi)|\,d\xi\\
&&\nonumber\qquad\leq\int_0^L\xi|g_1(t,\xi)-g_2(t,\xi)|\,d\xi\\
&&\label{1}\qquad+\left(\sign(u_1(t,0+)-u_2(t,0+))(u_1(t,0+)^2-u_2(t,0+)^2)/2\right)\\
&&\label{2}\qquad-\left(\sign(u^\infty_1(t)-u^\infty_2(t))(u^\infty_1(t)^2-u^\infty_2(t)^2)/2\right).
\end{eqnarray}
Note that the two last terms have the same absolute value. Without any loss of generality, we can assume $u_1(t,0+)\leq u_2(t,0+)$. If $u_2(t,0+)^2\leq u_1(t,0+)^2$,
then (\ref{1}) is negative, and the sum of (\ref{1}) and (\ref{2}) is non positive. Otherwise, $u_1(t,0+)^2\leq u_2(t,0+)^2$, and so $u_2(t,0+)>0$ and
$-u_1(t,0+)\leq u_2(t,0+)$. Then the BLN condition implies that $u_2^2/2=\int\xi g_2\,d\xi$, and so $F_2$ is a relaxation layer (thanks to Proposition \ref{theocouche}). This implies that $u_2^\infty=u_2$, and so
$u_2^\infty\geq u_1^\infty$. Hence (\ref{2}) is negative and the sum of (\ref{1}) and (\ref{2}) is non positive too.
This gives the result.\qquad\qed

We can now prove Theorem \ref{theoprolimite}.\vskip0.2cm

\noindent{\it Proof of Theorem \ref{theoprolimite}.} Let us first
show the stability of the system. We use Lemma \ref{lemm_coucheetbord} and Lemma
\ref{lemm_kineticaubord} with ($i$=1,2)
\begin{eqnarray*}
&&g_i(t,\xi\geq0)=f_i(t,0-,\xi\geq0),\qquad \xi\in(0,L), \ \ t\geq0,\\
&&g_i(t,\xi\leq0)=F_i(t,0+,\xi\leq0),\qquad \xi\in(-L,0), \ \ t\geq0.
\end{eqnarray*}
Adding the two estimates of those lemmas gives for any $t>0$
\begin{eqnarray*}
&&\qquad\int_{-\infty}^0\int_{-L}^L|f_1(t)-f_2(t)|\,dx\,d\xi+\int_0^\infty|u_1(t)-u_2(t)|\,dx
\\
&&\leq \int_{-\infty}^0\int_{-L}^L|f^0_1-f^0_2|\,d\xi\,dx+\int_0^{+\infty}|u^0_1-u^0_2|\,dx.
\end{eqnarray*}
This gives the stability and the uniqueness.

Let us show now the
existence. We fix a initial data $(f^0,u^0)$. We denote $H(T)$ the
set of boundary condition on the interface for $\xi<0$
$$
H=\{\xi g\in L^1(]0,T[\times]-L,0[), \ -1\leq g(t,\xi)\leq 0\}.
$$
Note that this set is convex.
We define the function $\mathcal{F}$ from $H$ to $H$ in the
following way. We consider the solution $f$ on
$(0,T)\times(-\infty,0)\times(-L,L)$ to the kinetic equation with
initial value $f^0$ and boundary condition $g$. Then we consider $u$
solution to the Burgers equation in $(0,T)\times(0,\infty)$ with
initial value $u^0$ and boundary condition (in the sense of BLN)
$\sqrt{2\int\xi f(t,0-,\xi)\,d\xi}$. Finally we consider $F$
solution to the layer problem with data $(f(t,0-,\xi),
(u(t,0+))^2/2)$. We then set
$$
\mathcal{F}(g)=\xi F(t,0+,\xi<0).
$$
Lemma \ref{lemm_kineticaubord} and
Lemma \ref{lemm_coucheetbord} ensure that $\mathcal{F}$ is
continuous from $H$ to $H$.

Consider a sequence of function $g_n\in H$ converging weakly to
$g\in H$. From Lemma \ref{lemm_kineticaubord}, we get that
$f_n(t,0-,\xi\geq0)$ converges strongly to $f(t,0-,\xi\geq0)$. Then,
 Lemma \ref{lemm_coucheetbord}
ensures that $\mathcal{F}(g_n)$ converges strongly in $H$ to
$\mathcal{F}(g)$. Hence $\mathcal{F}(H)$ is compact.

So, using the classical Schauder fixed point theorem, we get the
existence of a $g\in H$ such that $\mathcal{F}(g)=g$. The associated
functions $(f,F,u)$ is then solution to (\ref{limitefluide}), (\ref{limitecinetique}),
 (\ref{limitecouche}).\qed

\section{Asymptotic limit}

This section is devoted to the proof of Theorem
\ref{theolimitasymp}. Let us first show the following proposition.
It states that  the solution to (\ref{couplagecinetique}) for
$x\geq0$ converges to the the solution to (\ref{limitefluide}) with
the correct boundary condition.

\begin{prop}\label{theofluide}
Let  $f_\eps^0\in L^\infty(\R\times[-L,L])$ be such that
$0\leq\sign(\xi) f_\eps^0\leq1$. Denote  $\fe\in
L^\infty(\R^+\times\R\times[-L,L])$ the  solution to
(\ref{couplagecinetique}) with initial value $\fe^0$. Then there
exists  $\eps_n\to0$, two functions $f^0\in
L^\infty(\R\times[-L,L])$, $f\in L^\infty(\R^+\times\R\times[-L,L])$
such that  $\fe^0,\fe$ converge weakly to  $f^0,f$. Moreover the
function $u\in L^\infty(\R^+\times\R^+)$ defined by
\[
u(t,x)=\int_{-L}^{L}f(t,x,\xi)\,d\xi \ \ \mathrm{for} \ \
t\in\R^+,x\in\R^+,
\]
is solution to (\ref{limitefluide}).
\end{prop}
To show this Proposition, we first give a kinetic formulation of the $BLN$ condition (see \cite{Vovelle} for an other kind of kinetic formulation).
\begin{lemm}
Kinetic version of the BLN conditions. Consider
$g\in L^\infty(0,L)$ with $0\leq g(\xi)\leq 1$ and denote
$v=\sqrt{2\int_0^L\xi g(\xi)\,d\xi}$. We have
$$
u\stackrel{BLN}{:=}v
$$
if and only if there exists $h\in L^\infty(-L,0)$, $-1\leq
h(\xi)\leq 0$  and $m$ nonnegative measure on $(-L,L)$ such that
\begin{equation}\label{BLNcinetique}
 \xi(M(u,\xi)-(g+h)(\xi))= m'(\xi).
\end{equation}
\end{lemm}
\noindent{\bf Proof of the lemma}.
Indeed,  $u\stackrel{BLN}{:=}v$ if and only if   we have either
$u=v$ or $u\leq-v$. In the case $u=v$, we get (\ref{BLNcinetique})
with $h=0$ and $m(\xi)=\int_{-L}^\xi
\zeta(M(u,\zeta)-g(\zeta))\,d\zeta$. The function $m$ is equal to
zero for $\xi\leq 0$ and at $\xi=L$. It is increasing on $[0,u]$ and
decreasing on $[u,L]$. Hence it is nonnegative. If $u\leq-v$, we
take $h(\xi)=-\mathbf{1}_{[-\sqrt{u^2-v^2},0]}$ and
$m(\xi)=\int_{-L}^\xi \zeta(M(u,\zeta)-(g+h)(\zeta))\,d\zeta$. We
still have $m$ equal to 0 for $\xi\leq u$, $m(L)=0$ and $m$
increasing for $\xi<0$ and decreasing for $\xi>0$, so $m$ is still
nonnegative.

Conversely, assume that $u$ verifies (\ref{BLNcinetique}). Noting
that $\xi h\geq0$, and integrating in $\xi$, we find
$$
u^2/2\geq v^2/2.
$$
If $u\leq0$, this gives $u\leq -v$. Assume that $u\geq0$.
Multiplying (\ref{BLNcinetique}) by $(\xi)_-$ and integrating in
$\xi$, we find
$$
0\leq \int_{-L}^0\xi^2(-h)(\xi)\,d\xi=-\int_{-L}^0m(\xi)\,d\xi\leq0.
$$
Hence $h=0$, and integrating (\ref{BLNcinetique}) in $\xi$ we have
$$
u^2/2=v^2/2,
$$
and so $u=v$. This shows that if $u$ verifies (\ref{BLNcinetique})
then $u\stackrel{BLN}{:=}v$.\qquad\qed

\noindent{\bf Proof of the Proposition.} By weak compactness, there exists $\eps_n\to0$
such that $f_n=f_{\eps_n}$ converges weakly to $f$. Using averaging
lemmas, $\int f_n\,d\xi$ converges strongly to $\int f\,d\xi$, and
so $\M f_n$ converges strongly to $\M f$. Passing into the limit in
the kinetic domain gives that $f$ verifies for $x\in(-\infty,0)$
$$
\dt f+\xi\dx f=\M f-f.
$$
For the fluid domain, multiplying the equation by $\eps_n$ shows
that at the limit
\begin{equation}\label{mf=f}
\M f(t,x,\xi)-f(t,x,\xi)=0\qquad \mathrm{for \ almost \ every }\
t>0, \ x>0 \ \xi\in(-L,L).
\end{equation}
Thanks to Lemma \ref{simple}, there exists $m_n\geq0$ such that
$$
a_{\eps_n}(x)(\M f_n-f_n)=\partial_\xi m_n,
$$
and  those measures are uniformly bounded by
$$
\int_0^T\int_{-\infty}^\infty\int_{-L}^L
m_n(t,x,\xi)\,d\xi\,dx\,dt\leq \int_{-\infty}^\infty\int_{-L}^L\xi
f^0_{\eps_n}(x,\xi)\,d\xi\leq C.
$$
Thus, up to a subsequence, $m_n$ converges to nonnegative measure
$m$ (with possible concentration, especially at $x=0$). This gives
that for $x>0$:
$$
\dt f+\xi\dx f=\partial_\xi m, \qquad t>0, \ x>0,\ -L<\xi<L.
$$
This together with (\ref{mf=f}) is the kinetic formulation of the
Burgers equation and so $u(t,x)=\int f(t,x,\xi)\,d\xi$ verifies the
Burgers equation for $x>0$. We want now to recover the boundary
conditions on the Burgers equation. From the limit equation, we get
also that $f\in BV(\R_x;W^{-1,2}(\R^+_t\times(-L,L)))$ and so $f$
has a limit from the left and a limit from the right at each point
$x\in \R$ (in the sense of distribution). Moreover
$$
\xi(f(t,0+,\xi)-f(t,0-,\xi))=\partial_\xi \tilde{m}(t,0,\xi),\qquad t>0,
\xi\in (-L,L).
$$
The strong trace theorem  ensures that $f(t,0+,\xi)=M(u(t,0+),\xi))$
for $t>0$ and $-L<\xi< L$, and the kinetic formulation of the BLN
conditions (\ref{BLNcinetique}) gives
$$
u(t,0+)\stackrel{BLN}{:=}\sqrt{\int_0^L\xi f(t,0-,\xi)\,d\xi}.
$$
\qed

We now show the Theorem \ref{theolimitasymp} in the case of the
first set of initial values. We first consider the special initial
value, for $u_0>0$
$$
\psi^0(x,\xi)=\mathbf{1}_{\{x<0\}}M(u_0,\xi)+\mathbf{1}_{\{x>0\}}M(-u_0,\xi).
$$
We want to show that for this initial value, at the limit $\eps\to
0$ we get $f(t)=\psi^0$ for $t>0$. For this matter, we first
consider the initial value
$$
\psi^0_\delta(x,\xi)=\psi^0(x-\delta,\xi).
$$
Consider the limit function $f_\delta$. We have seen that
$u_\delta=\int f_\delta\,d\xi$ verifies the Burgers equation for
$x>0$ with initial value
$u^0_\delta(x)=\mathbf{1}_{\{x<\delta\}}u_0-\mathbf{1}_{\{x>\delta\}}u_0$.
Thanks to the finite speed of propagation there exists a finite time
(at least $\delta/u_0$) such that $u(t,0+)>0$. Hence, the BLN
conditions ensures that on this lapse of time
$u_\delta(t,0+)^2/2=\int_0^L\xi f_\delta(t,0-,\xi)\,d\xi$. But the
conservation of flux at the interface ensures then that
$f_\delta(t,0-,\xi<0)=0$ for this lapse of time and $\xi<0$. And so,
solving the equation on the kinetic domain shows that the solution
is stationary. The values at the interface does not change on time
either, and so the solution of the fluid domain is constant also.
Finally we find $f(t)=\psi^0_\delta$ for all $t<\delta/u_0$.
iterating the argument shows that this is true for all time $t>0$.
Using Lemma \ref{cadecroit}, we can show that for any $\eps$ and
$t>0$
$$
\int_{-\infty}^\infty\int_{-L}^L|f_\eps(t)-f^\delta_\eps(t)|\,d\xi\,dx\leq\int_{-\infty}^\infty\int_{-L}^L|\psi^0-\psi^0_\delta|\,d\xi\,dx=2\delta
u_0.
$$
up to subsequence, the functions converges strongly, and so
$$
\int_{-\infty}^\infty\int_{-L}^L|f(t)-\psi_\delta^0|\,d\xi\,dx\leq2\delta
u_0.
$$
Passing into the limit when $\delta\to0$ gives that $f(t)=\psi^0$
for every $t>0$.

Then, using Lemma \ref{cadecroit}, we can show that if $f_1^0\leq
f_2^0$ then for any $\eps$, $f_{1,\eps}\leq f_{2,\eps}$ (see \cite{Perthame}). And so, at
the limit $f_1\leq f_2$. The hypothesis on the initial data are
exactly $f^0\geq \psi^0$, hence, at the limit:
$$
f(t,x,\xi)\geq \psi^0(x,\xi)\qquad \mathrm{for} \ t>0, \ x\in\R, \
\xi\in(-L,L).
$$
This implies that $f(t,0-,\xi<0)=0$ for $t>0$, $\xi<0$. But the
condition implies also that the kinetic layer has to be of the
relaxation type. The condition is then the good one.

We are now left to show the asymptotic limit for the second set of
initial conditions. We need, in this case, to pass into the limit in
the term $F_\eps$. To avoid the possibility for the layer to be
swept away, we need to show that the shock is well confined against
the interface $x=0$. This will be provided by Proposition \ref{prop_confinement}. The hypothesis on the initial data
(\ref{initialchoc}), and a comparison argument gives that,
uniformly with respect to $\eps$
\[
|F_\eps(t,y,\xi)-M(-u^+;\xi)| \leq|F_g(y,\xi)-M(-u^+;\xi)| \
\mathrm{for} \ t\in\R^+,y\in\R^+,\xi\in\R,
\]
where  $F_g$ is the solution to the layer problem  (\ref{Milne})
associated to  $(g,V)$ with $g(\xi)=M(u^+-\eta;\xi)$
and $V=(u^+)^2/2$.

We need now to get strong compactness to pass into the limit in the
layer. For this, we use a blow-up technique first introduced in \cite{V2}.

\subsection{The ``blow up" method}\label{sectionlocal}

The equation on $F_\eps$ does not control the oscillations in time
when $\eps\to0$. The idea is to get back the balanced structure of
(\ref{BGK}) doing a zoom in time of the equation. We introduce a new
local variable of time $s\in\R$ and a new rescaled function defined
by
\begin{equation}\label{zoom}
\Fe(t,s,y,\xi)=F_\eps(t+\eps s,y,\xi).
\end{equation}
For almost every   fixed $t>0$, the function  $\Fe$ verifies
\begin{equation}\label{eqzoom}
\begin{array}{l}
\dis {\ds \Fe+\xi \dy\Fe=\M \Fe-\Fe, \ \ \mathrm{for} \ \
s\in]-t/\eps,+\infty[, y\in\R^+, \xi\in[-L,L],}\\[3mm]
\dis{\Fe(t,s,0,\xi)= \fe(t+\eps s,0,\xi) \ \ \mathrm{for} \ \
s\in]-t/\eps,+\infty[,  \xi\in[-L,L],}\\[3mm]
\dis{|\Fe(t,s,y,\xi)-M(-u^+;\xi)|\leq|F_g(y,\xi)-M(-u^+;\xi)|}\\[3mm]
\qquad\qquad\qquad \ \ \mathrm{for} \ \
s\in]-t/\eps,+\infty[,y\in\R^+,\xi\in[-L,L].
\end{array}
\end{equation}
Results obtained on the rescaled functions can be translated on the non-rescaled one (and vice versa) thanks to the following lemma whose proof can be found in \cite{V2}:
\begin{lemm}\label{local-global}(From local to global)
Let $F_{\eps_n}$, $F\in L^\infty(\R^p\times\R^q)$, $p,q$ integers. Then,
$F_{\eps_n}$ converges strongly to $F$ in
$L^1_{\mathrm{loc}}(\R^p\times\R^q)$ if and only if for
every $R_1>0$, $R_2>0$ and $R_3>0$:
\begin{equation}
\int_{B_p(0,R_1)}\int_{B_q(0,R_2)}\int_{B_p(0,R_3)}\left|F_{\eps_n}(x+\eps_n
y,z)-F(x,z)\right|\,dy\,dz\,dx
\stackrel{n\to+\infty}{\longrightarrow}0.
\end{equation}
\end{lemm}
Let us first apply this lemma on $f_\eps(\cdot,0,\cdot)$ with $p=1$ (time variable) and $q=1$ ($\xi$ variable). From Lemma \ref{lemm_kineticaubord}, the function $f_\eps(t,0,\xi>0)$ converges strongly to $f(t,0,\xi>0)$ in $L^1_{\mathrm{loc}}(\R^+\times(0,L))$. Hence for any $T>0$, $R>0$
$$
\int_0^T\left(\int_{0}^{R}\int_{0}^L|\Fe(t+\eps s,0,\xi)-f(t,0,\xi)|\,d\xi\,ds\right)\,dt\stackrel{\eps\to 0}{\longrightarrow}0.
$$
Hence, there is a subsequence $\eps_n$ such that for almost every fixed $t>0$,
the local function (restricted to positive $\xi$) $\Fen(t,\cdot,0,\cdot)$
converges strongly in $L^1_{\mathrm{loc}}(\R\times(0,L))$
to $f(t,0+,\cdot)$.

Let us fix such a  time $t>0$. We pass into the limit when $\eps_n$ goes
to 0.
For this we first show the following proposition.
\begin{prop}\label{prop_couche limite} Let $s_n$ be a sequence converging to $-\infty$ and consider
 $\Fn\in L^{\infty}((s_n,\infty)\times \R^+\times(-L,L))$, with $0\leq\sign(\xi)\Fn\leq1$  a sequence of solutions to
\begin{eqnarray}
\label{ajout_eq}&&\partial_s\Fn+\xi \dy\Fn=\M \Fn-\Fn,\\ %\qquad \mathrm{for} \ \ s\in (s_0,+\infty), y\in\R^+, \xi\in[-L,L],\\
\label{ajout_infini}&&|\Fn(s,y,\xi)-M(-u^+;\xi)|\leq|F_g(y,\xi)-M(-u^+;\xi)|\\
\nonumber&&\qquad\qquad\qquad\qquad \ \ \mathrm{for} \ \
s\in]s_n,+\infty[,y\in\R^+,\xi\in[-L,L].
\end{eqnarray}
Assume, in addition that $\Fn(\cdot,0,\cdot)$ (restricted to positive $\xi$) converges strongly in $L^1_{loc}(\R\times(0,L))$. Then, up to a subsequence, $\Fn$ converges strongly in $L^1_{\mathrm{loc}}(\R\times\R^+\times(-L,L))$ to a solution $F_l$ of (\ref{ajout_eq}) (\ref{ajout_infini}). In addition,  $\Fn(\cdot,0,\cdot)$ converges strongly in $L^1_{\mathrm{loc}}(\R\times(-L,L))$
to $F_l(\cdot,0,\cdot)$.
\end{prop}

\noindent{\bf Proof of the proposition.}
The function $\Fn$ is bounded by 1, so,
up to a subsequence, it converges weakly to a function denoted
$F_l$. We can now use the averaging lemmas on
the kinetic equation (\ref{ajout_eq}). We deduce that  $\M\Fn$ converges to $\M
F_l$ in $L^1_{\mathrm{loc}}$. Hence $F_l$ verifies the same equation (\ref{ajout_eq}).
Moreover, there exists a subsequence such that $M\Fn$ and  the restriction to $\xi>0$ of the trace $\Fn(\cdot,0,\cdot)$ converges almost everywhere.
For $\xi>0$, integrating the equation (\ref{ajout_eq}) along characteristics, we find
that for almost any $s,y$ and $n$ big enough:
$$
\Fn(s,y,\xi)=e^{-\frac{y}{\xi}}\Fn(s-\frac{y}{\xi},0,\xi)+\int_0^y\frac{e^{\frac{z-y}{\xi}}}{\xi}\M\Fn(s+\frac{z-y}{\xi},z,\xi)\,dz.
$$
From the  convergence of $\Fn$ almost everywhere  on the trace $(y=0,\xi>0)$ and the convergence of $\M \Fn$ almost everywhere, we deduce the convergence for almost every  $(s,y,\xi)\in \R\times\R^+\times(0,L)$ of $\Fn(s,y,\xi)$. This function is bounded by 1. The Lebesgue's Theorem ensures  the strong convergence of (the restriction for $\xi>0$ of) $\Fn$ in $L^1_{\mathrm{loc}}(\R\times\R^+\times(0,L))$.

For $\xi<0$, we use the same strategy, using (\ref{ajout_infini}).
Thanks to Lemma \ref{prop_proprieteslayer}, for any $\eta>0$, there exists $Y>0$ (as big as we wish) such that  for any $y>Y$ we have
$$
\int_{-L}^0|F_g(y,\xi)-M(-u^+;\xi)|\,d\xi\leq \eta.
$$
Thanks to (\ref{ajout_infini}), we deduce that for any $y>Y$,  we have for any $n,m$
$$
\int_{-L}^0 \sup_{s,s'}|\Fn(s,y,\xi)-\overline{F}_m(s',y,\xi)|\,d\xi\leq 2\eta.
$$
Integrating, again, the equation along the characteristics, we find for $n$ big enough
$$
\Fn (s,y,\xi)=e^{\frac{Y-y}{\xi}}\Fn(s+(Y-y)/\xi,Y,\xi)+\int_{Y}^y\frac{e^{(z-y)/\xi}}{\xi}\M \Fn (s+\frac{z-y}{\xi},z,\xi)\,dz.
$$
So, for $n,m$ big enough, we have for every $s$, and $y<Y$
\begin{eqnarray*}
&&\qquad\qquad\int_{-L}^0|\Fn (s,y,\xi)-\overline{F}_m(s,y,\xi)|\,d\xi\\
&&\leq 2\eta+\int_{-L}^0\int_{Y}^y\frac{e^{(z-y)/\xi}}{\xi}|\M \Fn-\M \overline{F}_m| (s+\frac{z-y}{\xi},z,\xi)\,dz\,d\xi.
\end{eqnarray*}
using the convergence of $\M\Fn$ almost everywhere, and the Lebesgue's Theorem, we find that for almost every $s\in \R$, and every  $y<Y$, the restriction for $\xi<0$ of $\Fn$ is Cauchy in $L^1_{\mathrm{loc}}(-L,0)$.
Using the Lebesgue's Theorem, and the previous result for $\xi>0$, we get the convergence of $\Fn$ in $L^1_{\mathrm{loc}}(\R\times\R^+\times(-L,L))$.

%Finally, $\Fn/\xi$ is bounded in the Lipshitz space in $y$ with values in distribution %$Lip(\R^+_y;W^{-1,\infty}_{\mathrm{loc}}(s_n,\infty)\times [(-L,0)\cup(0,L)])$.
%So by Aubin's lemma, we have the convergence on the trace: the limit of the trace is %the trace of the limit.

The result with $y=0$ gives  also the strong convergence of the trace.\qquad \qed
\vskip0.3cm
For almost every fixed $t>0$, we apply this proposition on $\Fen(t,\cdot,\cdot,\cdot)$. Passing
into the limit in (\ref{eqzoom}), we find that, up to a subsequence, $\Fen(t,\cdot,\cdot,\cdot)$ converges strongly to $\Fi(t,\cdot,\cdot,\cdot)$ solution to
\begin{equation}\label{eqzoomlimite}
\begin{array}{l}
\dis {\ds \Fi+\xi \dy\Fi=\M \Fi-\Fi, \ \ \mathrm{for} \ \
s\in\R, y\in\R^+, \xi\in[-L,L],}\\[3mm]
\dis{\Fi(t,s,0+,\xi)= f(t,0+,\xi) \ \ \mathrm{for} \ \
s\in\R,  \xi\in[0,L],}\\[3mm]
\dis{|\Fi(t,s,y,\xi)-M(-u^+;\xi)|\leq
|F_g(y,\xi)-M(-u^+;\xi)|,}\\[3mm]
\qquad\qquad\qquad\ \ \mathrm{for} \ \ s\in\R,y\in\R^+,\xi\in[-L,L].
\end{array}
\end{equation}
The following Liouville lemma (proven below) completely characterizes the limit.
\begin{lemm}\label{liouville}
There exists a unique solution  $\Fi(t,\cdot,\cdot,\cdot)$ to
(\ref{eqzoomlimite}). This solution does not depend on  $s$.
\end{lemm}
We  then deduce that  $\Fi$ is the unique solution  $F$ to the
kinetic layer problem (\ref{limitecouche}) associated to
$(f(t,0+,\xi\geq0),V=(u^+)^2/2)$. By the uniqueness of the limit,
the whole sequence is converging in $L^1_{\mathrm{loc}}$.

Lemma \ref{local-global} ensures the  convergence of the global
functions  $F_{\eps_n}$ from the result of convergence of the local
functions $\Fen$.

\subsection{Liouville's lemma}\label{sectionliouville}

This subsection is dedicated to the proof of Lemma \ref{liouville}.
Existence of a (steady) solution to (\ref{eqzoomlimite}) is given by Proposition \ref{theocouche} and a comparison principle. Assume that there exists two such solutions $F_1$ and $F_2$.
We first have:
\begin{equation}\label{zoomunicite}
\begin{array}{l}
\dis {\ds |F_1-F_2|+\xi \dy|F_1-F_2|=-\mathcal{D}(F_1,F_2), \
\mathrm{for} \
s\in\R, y\in\R^+, \xi\in\R,}\\
|F_1-F_2|(s,0+,\xi)= 0 \ \ \mathrm{for} \ \
s\in\R,  \xi\in\R^+,\\
|F_1-F_2|(s,y,\xi)\leq 2|F_g(y,\xi)-M(-u^+,\xi)| \ \ \mathrm{for} \
\ s\in\R, y\in\R^+, \xi\in\R,
\end{array}
\end{equation}
where
$$
\int_{-L}^L\mathcal{D}(F_1,F_2)\,d\xi=-\int_{-L}^L\sign(F_1-F_2)[\M F_1-\M F_2-F_1+F_2]\,d\xi\geq0,
$$
thanks to Lemma \ref{cadecroit}.
Thanks to Proposition
\ref{prop_confinement}, $F_g-M(-u^+;\cdot)$ is integrable, hence
$|F_1-F_2|(s,\cdot,\cdot)$ is integrable for every fixed $s$.
Integrating the first equation in  (\ref{zoomunicite}) with respect
to $y$ and $\xi$, we get:
\begin{equation}\label{rajout}
\begin{array}{l}
\qquad\dis{\frac{d}{ds}\int_{-L}^L\int_0^\infty|F_1-F_2|(s,y,\xi)\,d\xi\,dy}\\[3mm]
=\dis{-\int_{-L}^L\int_0^\infty\mathcal{D}(F_1,F_2)(s,y,\xi)\,d\xi\,dy+\int_{-\infty}^0\xi|F_1-F_2|(s,0+,\xi)\,d\xi\leq0.}
\end{array}
\end{equation}
We can deduce that
$\int_0^\infty\int_{-\infty}^{+\infty}|F_1-F_2|(s,y,\xi)\,d\xi\,dy$
is a non increasing bounded function. We denote   $\delta$ its limit
at $-\infty$. Consider now the functions translated in time
$F_i^n(s,\cdot,\cdot)=F_i(s-n,\cdot,\cdot)$ for $i=1,2$. Thanks to Proposition \ref{prop_couche limite}, up to a
subsequence, those functions converges in $L^1_\mathrm{loc}$ to
$F_1^\infty,F_2^\infty$ solutions to (\ref{eqzoomlimite}). In
particular, the quantity $|F_1^\infty-F_2^\infty|$ verifies
(\ref{zoomunicite}). Thanks to the Lebesgue's dominated convergence
theorem, we get
$\int_0^\infty\int_{-\infty}^{+\infty}|F^\infty_1-F^\infty_2|(s,y,\xi)\,d\xi\,dy=\delta$
for every $s\in\R$. Hence, we have
$$
\int_0^\infty\int_{-L}^L \mathcal{D}(F_1^\infty,F^\infty_2)\,d\xi\,dy=0,
$$
and
$$
|F^\infty_1-F^\infty_2|(s,0,\xi)=0.
$$
The proof now follows the lines of the proof of Lemma \ref{lemm_uniqueness}.
Thanks to Lemma \ref{cadecroit}, we have
$\partial_\xi(\sign(F^\infty_1-F^\infty_2))=0$. Let us fix $s\in\R$, $y\in\R^+$ and $\xi<0$.
We denote
$$
\tilde{F}_i(\tau)=F^\infty_i(s+\tau,y+\tau \xi,\xi),\qquad i=1,2, \ \ \ \tau\in \R,
$$
and
$$
\Omega_{s,y,\xi}=\{0\leq\tau\leq -y/\xi \setminus \tF_1(\tau)>\tF_2(\tau)\}.
$$
Note that this set is open  since $\tF_1-\tF_2$ is Lipschitz (and so continuous). Assume that it is not empty. Denote $\tau$ one of its elements and $\tau_0=\sup\{\lambda>\tau \setminus (\tau,\lambda)\subset \Omega_{s,y,\xi}\}$.
For every $\lambda\in(\tau,\tau_0)$ and every
$\zeta\in(-L,L)$ we have: $F^\infty_1(s+\lambda,y+\lambda\xi,\xi)> F^\infty_2(s+\lambda,y+\lambda\xi,\xi)$, hence we have
$F^\infty_1(s+\lambda,y+\lambda\xi,\zeta)\geq F^\infty_2(s+\lambda,y+\lambda\xi,\zeta)$ too. This leads to:
\[
\partial_\tau (\tF_1-\tF_2)+(\tF_1-\tF_2)\geq \M \tF_1-\M \tF_2\geq0, \qquad
\mathrm{on \ } (\tau,\tau_0).
\]
Hence, we find that $\tau_0=-y/\xi$ and $(F^\infty_1-F^\infty_2)(s-y/\xi,0,\xi)\geq
(\tF_1-\tF_2)(\tau)e^{y/\xi}>0 $. This gives a contradiction. Hence
$\Omega_\xi$ is empty and so $F^\infty_1\leq F^\infty_2$ for $\xi<0$. Exchanging
the indices gives:
\[
F^\infty_1(s,y,\xi)=F^\infty_2(s,y,\xi) \qquad s\in\R, y\in\R^+, \ \xi\in(-L,0).
\]
%For $\xi>0$ we consider
%$$
%\Omega_{s,y,\xi}=\{-y/\xi\leq\tau\leq0 \setminus \tF_1(\tau)<\tF_2(\tau)\}.
%$$
%In the same way, we show that, if $\Omega_{s,y,\xi}$ is not empty, then $(-y/\xi,0)$ is
%included in $\Omega_{s,y,\xi}$ and
%$(F_1-F_2)(s-y/\xi,0,\xi)\leq
%(\tF_1-\tF_2)(0)e^{y/\xi}>0 $.

Hence, for every $y>0$
\begin{equation}\label{rerajout}
\begin{array}{l}
\qquad\dis{\frac{d}{ds}\int_{-L}^L\int_0^y|F^\infty_1-F^\infty_2|(s,z,\xi)\,d\xi\,dz}\\[3mm]
=\dis{-\int_{0}^L\xi|F^\infty_1-F^\infty_2|(s,y,\xi)\,d\xi=D(s,y)\leq0.}
\end{array}
\end{equation}
So for any $y>0$, $\int_0^y\int_{-L}^L|F^\infty_1-F^\infty_2|(s,z,\xi)\,d\xi\,dz$ is
a  non increasing bounded function. We denote $d(y)$ its limit for $s$ goes to $-\infty$. Consider a new translation in time
$$
F_i^{\infty,n}(s,\cdot,\cdot)=F_i^\infty(s-n,\cdot,\cdot),
$$
for $i=1,2$. As before, up to a subsequence, those functions converge in $L^1_{\mathrm{loc}}$ to $F_1^{\infty,\infty},F_2^{\infty,\infty}$. They verify
$$
|F_1^{\infty,\infty}-F_2^{\infty,\infty}|(s,y,\xi<0)=0,
$$
and for any $y>0$, and $s\in\R$ $\int_{-L}^L\int_0^y|F^{\infty,\infty}_1-F^{\infty,\infty}_2|(s,z,\xi)\,d\xi\,dz=d(y)$.
Hence, from (\ref{rerajout}), for any $y>0$, $s\in \R$:
$$
\int_{0}^L\xi|F^{\infty,\infty}_1-F^{\infty,\infty}_2|(s,y,\xi)\,d\xi=0.
$$
Finally $F^{\infty,\infty}_1=F^{\infty,\infty}_2$ and $d(y)=0$ for any $y$.
 We
can deduce that  $|F^\infty_1-F^\infty_2|(s,y,\xi)=0$ everywhere,
and so $\delta=0$. This gives  $F_1=F_2$.\qed

\vskip0.2cm
\noindent{\bf Acknowledgements:} This work was partially supported by the NSF Grant DMS:0607053.\\
We thank the referee for providing a much more elegant proof of Proposition \ref{prop_confinement}.

\bibliography{biblioprague}

\end{document}